\documentclass[a4paper,12pt]{article}

\usepackage{amsmath,amsthm,amssymb,mathrsfs,latexsym,amsfonts}
\usepackage{graphicx,psfrag,epsfig}
\usepackage{bbm}
\usepackage{a4wide}
\usepackage[english]{babel}
\usepackage[latin1]{inputenc}
\usepackage{authblk}

\numberwithin{equation}{section}

\newtheorem{theorem}{Theorem}[section]

\newtheorem{proposition}[theorem]{Proposition}

\newtheorem{Main}{Theorem}

\newcounter{paraga}[section]

\newcommand{\N}{\mathbb{N}}
\newcommand{\Z}{\mathbb{Z}}
\newcommand{\Q}{\mathbb{Q}}
\newcommand{\R}{\mathbb{R}}
\newcommand{\C}{\mathbb{C}}

\begin{document}

\def\MP{\,{<\hspace{-.5em}\cdot}\,}
\def\SP{\,{>\hspace{-.3em}\cdot}\,}
\def\PM{\,{\cdot\hspace{-.3em}<}\,}
\def\PS{\,{\cdot\hspace{-.3em}>}\,}
\def\EP{\,{=\hspace{-.2em}\cdot}\,}
\def\PP{\,{+\hspace{-.1em}\cdot}\,}
\def\PE{\,{\cdot\hspace{-.2em}=}\,}
\def\N{\mathbb N}
\def\C{\mathbb C}
\def\Q{\mathbb Q}
\def\R{\mathbb R}
\def\T{\mathbb T}
\def\A{\mathbb A}
\def\Z{\mathbb Z}
\def\demi{\frac{1}{2}}

\begin{titlepage}
  \title{\LARGE{\textbf{Some remarks on the Classical KAM Theorem, following P{\"o}schel}}}
  \author{Abed Bounemoura \\
    CNRS - PSL Research University\\
    (Universit{\'e} Paris-Dauphine and Observatoire de Paris)}
\end{titlepage}

\maketitle

\begin{abstract}
We propose a slight correction and a slight improvement on the main result contained in ``A lecture on Classical KAM Theorem" by J. P{\"o}schel. 
\end{abstract}

\section{Introduction and results}\label{s1}

The paper~\cite{Pos01} contains a very nice exposition of the classical KAM theorem which has been very influential. It is our purpose in this short and non self-contained note to add two remarks to this remarkable paper. 

The first one concerns a technical mistake\footnote{The choices of $h_0$ and $K_0$, page 23 in~\cite{Pos01}, violate the condition $h_0 \leq \alpha(2K_0^\nu)^{-1}$.} in the proof of the main abstract statement Theorem A, which has been recently pointed out and corrected in the PhD thesis~\cite{Kou19}. Yet a correction of this mistake, following P{\"o}schel arguments, leads to a final statement which is both less elegant and quantitatively weaker. We would like to explain how, by modifying slightly the arguments using ideas due to R{\"u}ssmann (see for instance~\cite{Rus01}), Theorem A of~\cite{Pos01} can be proved without any changes. The aforementioned modifications consist of replacing the crude Fourier truncation by a more refined polynomial approximation, and then set an iterative scheme with a linear\footnote{We would like to quote here the paper~\cite{Rus89}: ``It has often been said that the rapid convergence of the Newton iteration is necessary for compensating the influence of small divisors. But a deeper analysis shows that this is not true. The Newton method compensates not only the influence of small divisors but also many bad estimates veiling the true structure of the problems."}, rather than super-linear, speed of convergence.

The second one concerns the application of Theorem A to an $\varepsilon$-perturbation of a non-degenerate integrable Hamiltonian system. This gives persistence of a set of positive measure of analytic invariant quasi-periodic tori with fixed diophantine frequencies, such that each torus in this set is at a distance of order $\sqrt{\varepsilon}$ to its associated unperturbed invariant torus. By using a more adapted version of Theorem A, we can actually show that the distance is of order $\varepsilon/\alpha$, where $\alpha$ is the constant of the Diophantine vector. This is not a new result, as this was already proved in~\cite{Vil08} using a refinement of Kolmogorov approach (for an individual torus).

\medskip

So let us recall the main result of~\cite{Pos01}, keeping the same notations. For a given domain $\Omega \subseteq \R^n$, consider a subset $\Omega_\alpha \subseteq \Omega$ of Diophantine vectors with constant $\alpha>0$ and exponent $\tau \geq n-1$. Given $0<r,s,h\leq 1$, define
\[ D_{r,s}=\{I \; | \; |I|<r\} \times \{\theta \; | \; |\mathrm{Im}(\theta)|<s\} \subseteq \C^n \times \C^n, \quad O_h=\{\omega \; | \; |\omega-\Omega_\alpha|<h\} \subseteq \C^n \]
where $|\,.\,|$ is the sup norm for vectors, and let $|\,.\,|_{r,s,h}$ the sup norm for functions defined on $D_{r,s} \times O_h$ and $|\,.\,|_{L}$ the Lipschitz semi-norm with respect to $\omega$. Let $N(I,\omega)=e(\omega)+\omega\cdot I$, which can be seen as a family $N_\omega$ of linear integrable Hamiltonian depending on parameters $\omega \in \Omega$; the family of embedding $\Phi_0 : \T^n \times \Omega \rightarrow \R^n \times \T^n$ defined by $\Phi_0(\theta,\omega)=(0,\theta)$ defines, for each $\omega \in \Omega$, a Lagrangian torus invariant by the Hamiltonian flow of $N_\omega$ and quasi-periodic of frequency $\omega$.

\begin{Main}\label{th1}
Let $H=N+P$. Suppose $P$ is real-analytic on $D_{r,s} \times O_h$ with 
\begin{equation}\label{seuil1}
|P|_{r,s,h} \leq \gamma \alpha r s^\nu, \quad \alpha s^\nu \leq h
\end{equation} 
where $\nu=\tau+1$ and $\gamma$ is a small constant depending only on $n$ and $\tau$. Then there exist a Lipschitz map $\varphi : \Omega_\alpha \rightarrow \Omega$ and a Lipschitz familiy of real-analytic Lagrangian embedding $\Phi : \T^n \times \Omega_\alpha \rightarrow \R^n \times \T^n$ that defines, for each $\omega \in \Omega_\alpha$, a Lagrangian torus invariant by the Hamiltonian flow of $H_{\varphi(\omega)}$ and quasi-periodic of frequency $\omega$. Moreover, $\Phi$ is real-analytic on $T_*= \{\theta \; | \; |\mathrm{Im}(\theta)|<s/2\}$ for each $\omega$ and
\begin{equation}\label{est1}
\begin{cases}
|W(\Phi-\Phi_0)|, \, \alpha s^\nu |W(\Phi-\Phi_0)|_L \leq c(\alpha r s^\nu)^{-1}|P|_{r,s,h} \\
|\varphi-\mathrm{Id}|, \, \alpha s^\nu |\varphi-\mathrm{Id}|_L \leq cr^{-1}|P|_{r,s,h}
\end{cases}
\end{equation}
uniformly on $T_* \times \Omega_\alpha$ and $\Omega_\alpha$ respectively, where $c$ is a large constant depending only on $n$ and $\tau$, and $W=\mathrm{Diag}(r^{-1}\mathrm{Id},s^{-1}\mathrm{Id})$.
\end{Main}

As expressed in~\eqref{est1}, the map $(\Phi,\varphi)$ is Lipschitz regular with respect to $\omega \in \Omega_\alpha$, and its Lipschitz norm (suitably weighted) is close to the one of $(\Phi_0,\mathrm{Id})$; this is all what is needed to transfer the positive measure in parameter space $\omega \in \Omega_\alpha$ to a positive measure of quasi-periodic solutions in phase space. One course one may ask whether $(\Phi,\varphi)$ is more regular with respect to $\omega \in \Omega_\alpha$ (since $\Omega_\alpha$ is a closed set, smoothness has to be understood in the sense of Whitney). In fact, the sketch of proof we will give below implies the following: given any $l \in [1,+\infty[$, provided~\eqref{seuil1} is replaced by
\[ |P|_{r,s,h} \leq \gamma(l) \alpha r s^\nu \]     
for some $h>0$ and some $\gamma(l)>0$, $(\Phi,\varphi)$ is of class $C^l$ with respect to $\omega$: we simply chose $l=1$ in Theorem~\ref{th1} to obtain Lipschitz regularity. However, as $l \rightarrow +\infty$, $\gamma(l) \rightarrow 0$ and thus we cannot conclude that $(\Phi,\varphi)$ is smooth. In order to reach such a statement, one can replace the linear scheme of convergence by the usual super-linear scheme (as described in~\cite{Pos01} for instance) but then the exponent $\nu$ in~\eqref{seuil1} has to be deteriorate: given any $\mu>\nu$, we have that $(\Phi,\varphi)$ is smooth with respect to $\omega$ provided~\eqref{seuil1} is replaced by
\[ |P|_{r,s,h} \leq \gamma(\mu,\nu) \alpha r s^\mu \]     
for some $h>0$ and some $\gamma(\mu,\nu)>0$: again $\gamma(\mu,\nu) \rightarrow 0$ as $\mu \rightarrow \nu$. Popov (see~\cite{Pop04}) showed that one can even go further and obtain some Gevrey smoothness of $(\Phi,\varphi)$ under a stronger smallness condition; without going into these rather technical issues, let us just say that $(\Phi,\varphi)$ can be shown to be Gevrey with exponent $1+\mu$ provided the polynomially small threshold $s^\nu$ in~\eqref{est1} is replace by a super-exponentially small threshold of order $\exp(-c(1/s)^a)$ with $a=a(\mu,\nu)=\nu/(\mu-\nu)$. This is probably the best smoothness one can achieve in general. 

\medskip

Next we consider a small perturbation of a non-degenerate integrable Hamiltonian, that is a real-analytic Hamiltonian of the form
\[ H(q,p)=h(p)+ f(q,p), \quad |f|\leq \varepsilon \] 
where $|f|$ is the sup norm on a proper complex domain. Introducing frequencies as independent parameters as in~\cite{Pos01}, one can write $H$ as in Theorem~\ref{th1} with
\[ P=P_f+P_h, \quad |P_f|\leq \varepsilon, \quad |P_h|\leq Mr^2 \]
where $M$ is a bound on the Hessian of $h$. At that point, the best choice for $r$ seems to be $r\simeq \sqrt{\varepsilon}$ so that the size of $P$ is of order $\varepsilon$ and Theorem~\ref{th1} can be applied; yet with such a choice it is obvious that because of the estimates for $\varphi$ in~\eqref{est1}, the distance between the perturbed and unperturbed torus will be of order $\varepsilon/r\simeq \sqrt{\varepsilon}$. Such an argument, used in~\cite{Pos01}, do not take into account the fact that the term $P_h$ is actually integrable and at least quadratic in $I$ (that is, $P_h(0,\omega)=0$ and $\nabla_I P_h(0,\omega)=0$): this is an important point, as the size of $P_h$ will effectively enter into the conditions~\eqref{seuil1} but not in the estimates~\eqref{est1}, simply because $P_h$ do not get involved in the approximation procedure nor contribute to the linearized equations one need to solve at each step of the iteration. Then, taking into account the estimate for $P_h$ (which itself is a consequence of the fact that it is at least quadratic in $I$), the requirement
\[ |P| \lesssim \alpha r s^\nu  \]
is then obviously implied by the conditions
\[ |P_f| \lesssim \alpha r s^\nu, \quad r \lesssim \alpha s^\nu \]
and thus we can state the following theorem (with a change of notations).

\begin{Main}\label{th2}
Let $H=N+P+Q$. Suppose $P$, $Q$ are real-analytic on $D_{r,s} \times O_h$, $Q$ is integrable and at least quadratic in $I$ with $|Q|_{r,h} \leq Mr^2$ and 
\begin{equation}\label{seuil2}
|P|_{r,s,h} \leq \gamma \alpha r s^\nu, \quad r \leq \delta M^{-1} \alpha s^\nu, \quad \alpha s^\nu \leq h
\end{equation} 
where $\nu=\tau+1$, $\gamma$ and $\delta$ are small constants depending only on $n$ and $\tau$. Then there exist a Lipschitz map $\varphi : \Omega_\alpha \rightarrow \Omega$ and a Lipschitz familiy of real-analytic Lagrangian embedding $\Phi : \T^n \times \Omega_\alpha \rightarrow \R^n \times \T^n$ that defines, for each $\omega \in \Omega_\alpha$, a Lagrangian torus invariant by the Hamiltonian flow of $H_{\varphi(\omega)}$ and quasi-periodic of frequency $\omega$. Moreover, the estimates~\eqref{est1} holds true.
\end{Main}

We may now choose $r$ as large as possible, namely $r \simeq \alpha s^{\nu}$, and obtain as a consequence that the distance between perturbed and unperturbed torus is of order $\varepsilon(\alpha s^{\nu})^{-1}$. As we already said, this fact was proved in~\cite{Vil08}; alternatively, a slight modification in the proof in~\cite{BF17} yields the same result.

\section{Sketch of proof}

In this section, we will sketch the proof of Theorem~\ref{th1} and Theorem~\ref{th2}; actually, we will simply indicate the modifications with respect to~\cite{Pos01} and we will use the same convention for implicit constants depending only on $n$ and $\tau$.

\begin{proposition}\label{step}
Let $H=N+P$, and suppose that $|P|_{s,r,h} \leq \varepsilon$ with
\begin{equation}\label{seuil}
\begin{cases}
\varepsilon \PM \alpha\eta^2 r \sigma^{\nu}, \\
\varepsilon \PM h r, \\
h \leq \alpha(2K^{\nu})^{-1}, \quad K \EP \sigma^{-1}\log(n\eta^{-2})
\end{cases}
\end{equation}
where $0<\eta<1/8$ and $0<\sigma<s/5$. Then there exists a real-analytic transformation 
\[ \mathcal{F}=(\Phi,\varphi) : D_{\eta r,s-5\sigma} \times O_{h/4} \rightarrow   D_{r,s} \times O_{h} \] 
such that $H \circ \mathcal{F}=N_++P^+$ with
\begin{equation}\label{estim1}
|P_+| \leq 9\eta^2\varepsilon
\end{equation}
and
\begin{equation}\label{estim2}
\begin{cases}
|W(\Phi-\mathrm{Id})|, \,  |W(D\Phi-\mathrm{Id})W^{-1}|\MP (\alpha r \sigma^\nu)^{-1}\varepsilon \\
|\phi-\mathrm{Id}|, \, h|D\varphi-\mathrm{Id}|_L \MP r^{-1}\varepsilon
\end{cases}
\end{equation}
uniformly on  $D_{\eta r,s-5\sigma} \times O_{h}$ and  $O_{h/4}$, with $W=\mathrm{Diag}(r^{-1}\mathrm{Id},\sigma^{-1}\mathrm{Id})$.
\end{proposition}

The above proposition is a variant of the KAM step of~\cite{Pos01}, which we already used in~\cite{Bou19}. The only difference is that in~\cite{Pos01}, instead of~\eqref{seuil} the following conditions are imposed 
\begin{equation}\label{seuilP}
\begin{cases}
\varepsilon \PM \alpha\eta r \sigma^{\nu}, \\
\varepsilon \PM h r, \\
h \leq \alpha(2K^{\nu})^{-1}
\end{cases}
\end{equation}
with a free parameter $K \in \N^*$, leading to the following estimate
\begin{equation}\label{estimP}
|P_+| \MP (\varepsilon(r\sigma^\nu)^{-1}+\eta^2+K^ne^{-K\sigma})\varepsilon.
\end{equation}
instead of~\eqref{estim1}. The last two terms in the estimate~\eqref{estimP} comes from the approximation of $P$ by a Hamiltonian $R$ which is affine in $I$ and a trigonometric polynomial in $\theta$ of degree $K$; to obtain such an approximation, in~\cite{Pos01} the author simply truncates the Taylor expansion in $I$ and the Fourier expansion in $\theta$ to obtain the following approximation error
\[ |P-R|_{s-\sigma,2\eta r, h} \MP (\eta^2+K^ne^{-K\sigma}).  \]
Yet we can use a more refined approximation result, which allows to get rid of the factor $K^n$ in the above estimate. More precisely, we use Theorem $7.2$ of~\cite{Rus01} (choosing, in the latter reference, $\beta_1=\cdots=\beta_n=1/2$ and $\delta^{1/2}=2\eta$ for $\delta \leq 1/4$); with the choice\footnote{There is a constant depending only on $n$ that we left implicit in the definition of $K$, which depends on the precise choice of norms for real and integer vectors,  see~\cite{Rus05} for instance.} of $K$ as in~\eqref{seuil}, this gives another approximation $\tilde{R}$ (which is nothing but a weighted truncation, both in the Taylor and Fourier series, which is affine in $I$ and of degree bounded by $K$ in $\theta$) and a simpler error
\[ |P-\tilde{R}|_{s-\sigma,2\eta r, h} \leq 8\eta^2.  \]
As for the first term in the estimate~\eqref{estimP}, it can be easily bounded by $\eta^2\varepsilon$ in view of the first part of~\eqref{seuil} which is stronger than the first part of~\eqref{seuilP} required in~\cite{Pos01}. 

Now, at variance with~\cite{Pos01}, we will use Proposition~\ref{step} in an iterative scheme with a linear speed of convergence as $\eta$ will be chosen to be a small but fixed constant: for convenience, let us set
\[ \eta=10^{-1}4^{-\nu}, \quad \kappa=9\eta^2. \]
Next, we define for $i \in \N$,
\[ \sigma_0=s/20, \quad \sigma_i=2^{-i}\sigma_0, \quad s_0=s, \quad s_{i+1}=s_i-5\sigma_i \]
so that $s_i$ converges to $s/2$. Then, for $K_i\EP \sigma_i^{-1}\log(n\eta^2)\EP \sigma_i^{-1}$, we set
\[ h_i=\alpha(2K_i^\nu)^{-1}=2^{-i \nu}h_0, \quad h_i\PE \alpha\sigma_i^{\nu} \]
and the condition $\alpha s^\nu \leq h$ implies in particular than $h_0\leq h$. Finally, we put
\[ \varepsilon_i=\kappa^i\varepsilon, \quad r_i=\eta^i r \]
and we verify that Proposition~\ref{step} can be applied infinitely many times: the third condition of~\eqref{seuil} holds true by definition, whereas the first two conditions of~\eqref{seuil} amount to $\varepsilon_i \PM \alpha r_i \sigma_i^\nu$ which, in view of our choice of $\eta$, holds true for all $i \in \N$ provided it holds true for $i=0$; for $i=0$ the condition is satisfied in view of the threshold $\varepsilon \leq \gamma \alpha r s^\nu$. Once we can iterate Proposition~\ref{step} infinitely many times, the convergence proof and the final estimates follow exactly as in~\cite{Pos01}, since the sequences $ \varepsilon_i(h_ir_i)^{-1}$and $\varepsilon_i(h_i^2r_i)^{-1}$ decrease geometrically, again by our choice of $\eta$. This concludes the sketch of proof.

To prove Theorem~\ref{th2}, one needs the following variant of Proposition~\ref{step}.

\begin{proposition}\label{step2}
Let $H=N+P+Q$, suppose that $|P|_{s,r,h} \leq \varepsilon$, $|Q|_{r,h}\leq M r^2$ with $Q$ integrable and at least quadratic in $I$ and
\begin{equation}\label{seuilstep2}
\begin{cases}
\varepsilon \PM \alpha\eta^2 r \sigma^{\nu}, \\
r \PM M^{-1} \alpha\eta^2 \sigma^{\nu}, \\
\varepsilon \PM h r, \\
h \leq \alpha(2K^{\nu})^{-1}, \quad K=n\sigma^{-1}\log(\eta^{-2})
\end{cases}
\end{equation}
where $0<\eta<1/4$ and $0<\sigma<s/5$. Then there exists a real-analytic transformation 
\[ \mathcal{F}=(\Phi,\varphi) : D_{\eta r,s-5\sigma} \times O_{h/4} \rightarrow   D_{r,s} \times O_{h} \] 
such that  $H \circ \mathcal{F}=N_++P_++Q$ with the estimates~\eqref{estim1} and~\eqref{estim2}.
\end{proposition}

Let $\tilde{R}$ be the approximation of $P$; if $\{.,.\}$ denotes the Poisson bracket and $[\,.\,]$ averaging over the angles, we solve the equation
\[ \{F,N\}=\tilde{R}+Q-[\tilde{R}+Q] \]
which, since $Q$ is integrable, is exactly the equation
\[ \{F,N\}=\tilde{R}-[\tilde{R}] \]
that is solved in~\cite{Pos01} (with, of course, $R$ instead of $\tilde{R}$ as we explained above). This justifies that the transformation in Proposition~\ref{step2} is the same as in Proposition~\ref{step}, and in particular it satisfy the estimates~\eqref{estim1}. The only difference is that the new Hamiltonian writes
\[ H \circ \mathcal{F}=N_++P_++Q, \quad N_+=N+[\tilde{R}] \]
with
\[ P_+=\int_{0}^1\{(1-t)[\tilde{R}]+t\tilde{R}+Q,F\} \circ X_F^t dt+(P-\tilde{R})\circ X_F^1. \]
As compared to~\cite{Pos01}, there is an extra term in $P_+$ coming from $Q$, whose contribution is easily bounded by the simple Poisson bracket
\[ |\{Q,F\}| \MP Mr(\alpha \sigma^\nu)^{-1}\varepsilon \]
and, in view of the extra condition we imposed in~\eqref{seuilstep2}, one can easily arrange the estimate~\eqref{estim2}. This justifies Proposition~\ref{step2}, and the iteration leading to Theorem~\ref{th2} is exactly the same as the one leading to Theorem~\ref{th1}. 

\bigskip

\textit{Acknowledgements.} The author have benefited from partial funding from the ANR project Beyond KAM.

\addcontentsline{toc}{section}{References}
\bibliographystyle{amsalpha}
\bibliography{RemarksKAMrevised}

\end{document}